\documentclass[12pt]{article}
\usepackage[psamsfonts]{amssymb}
\usepackage{latexsym}
\usepackage{amsmath,amssymb,amscd}
\newcommand\qed{{\unskip\nobreak\hfil\penalty50\hskip2em\vadjust{}
     \nobreak\hfil$\Box$\parfillskip=0pt\finalhyphendemerits=0\par}}

\newcommand{\x}{\times}
\newcommand{\ox}{\otimes}

\renewcommand{\>}{\rangle}

\renewcommand{\a}{\alpha}
\renewcommand{\b}{\beta}
\renewcommand{\d}{\delta}
\newcommand{\D}{\Delta}
\newcommand{\e}{\varepsilon}
\newcommand{\g}{\gamma}

\renewcommand{\l}{\lambda}
\renewcommand{\L}{\Lambda}
\newcommand{\n}{\nabla}
\newcommand{\var}{\varphi}
\newcommand{\s}{\sigma}

\renewcommand{\O}{\Omega}

\renewcommand{\i}{\infty}
\newcommand{\p}{\partial}

\title{Free Analysis Questions I:\\
Duality Transform for the Coalgebra of $\p_{X:B}$}
\author{Dan Voiculescu
\\ Department of Mathematics\\
University of California at Berkeley\\
Berkeley, CA\ \ 94720-3840\\
E-mail: dvv@math.berkeley.edu}
\date{June 6, 2003  \\ Preliminary Version}

\begin{document}
\maketitle

\begin{abstract} 
A duality transform for the coalgebra of the free difference quotient derivation-multiplication of an operator
with respect to a free algebra of scalars is constructed. The dual object is realized in an algebra  of
matricial analytic functions endowed with yet another generalization of the difference quotient derivation.

\end{abstract}

\setcounter{section}{0}
\section{Introduction}
\label{sec1}
The analysis aspects of variables with the highest degree of noncommutativity pose many dificult problems. 
The free difference quotient $\p_{X:B}$ is a derivation-comultiplication on the noncommutative polynomials
$B\<X\>$ with algebra of scalars $B$. In free probability, in questions related to random matrices and to the
 free analogue of entropy, $\p_{X:B}$ is the natural basic differential operator and like for the usual 
${\p}\slash{\p x}$, the questions are about the action on more general functions than polynomials 
(actually on noncommutatice $L^2$-functions).

In fact $\p_{X:B}$ is a coassociative comultiplication and in [12] we showed that this coalgebra and its dual,
are the key to the analytic subordination results in free probability([2], [10]). We showed that under some
technical  restrictions, the corepresentations, i.e. group-like elements, are generalized $B$-resolvents of
$X$.  Such generalized resolvents occur in the operator-valued extensions of free probability [9] and  have been
an essential ingredient in some remarkable random matrix work ([4], [7]). They also play an important role in
system theory (see [5] and references therein). Therefore, one side of the study of duality for the
coalgebra of
 $\p_{X:B}$, is  a generalization of spectral analysis, via resolvents, to the case of "free scalars``, i.e.
the  "scalars`` are an algebra of operators with the highest degree of noncommutativity w.r.t.the analyzed
operator.

In [12], we also took a look at the coalgebra for which the comultiplication is a derivation (GDQ---generalized 
difference quotient rings) and found that algebraically, the class has the remarkable property of being
self-dual, i.e. the dual comultiplication is again a derivation. In our functional analysis context, in the
simplest case of a  selfadjoint operator and scalars $ B= {\mathbb C}$ (i.e. commuting) the dual identifies
with an algebra of analytic functions on the resolvent and comultiplication is again given by the difference
quotient.

More recently we learned that the coalgebras with derivation-comultiplication had been noted by G.C. Rota(see
[6])
 based on the example of the difference quotient on commutative polynomials (without apparently realizing the
class  was closed under duality). The dual algebra in this combinatorial context is related to sequences of
classical
 polynomials and umbral calculus (see [1] for a recent free probability extension).

The present paper can be viewed as a continuation of [12].
 
On the algebraic side, we show that the multivariable situation $B\<X_1,\dots , X_n\>$ with several noncommuting 
variables and corresponding partial free difference quotients $\p_1,\dots ,\p_n$ can be reduced to
a one-variable  GDQ ring and can therefore be studied via coalgebra duality. We also improve our results on
corepresentations by identifying new ones. They are however, natural extensions of the generalized resolvents
we already found in [12]. We also show that free difference quotients are typical of GDQ rings because a
variable $X$, so that 
$\p X = 1\ox 1$, and the "scalars`` $N= {\rm Ker}\p$ always generate free polynomials on which $\p$ acts as
$\p_{X:N}$.

The main aim of this paper is to construct a suitable framework for the dual GDQ structure in the case of an 
operator $Y$ and a noncommutative algebra of scalars $B$. Approaching duality via a map of the dual $E'$ of the
Banach algebra containing $B$ and $Y$ into matrices indexed by corepresentations, we need a certain GDQ structure on
the matricial functions. Since in case $B={\mathbb C}$ the dual is a GDQ ring of analytic functions with respect
to the difference quotient on the resolvent set of $Y$, to deal with general $B$ requires a generalization of this.
It turns out that we need to consider collections of matricial objects at all levels, very much like in K-theory
or in the theory of operator spaces. Thus, for instance, instead of the scalar resolvent set, we will have an object
combining all matricial $B$-resolvent sets, tied together by natural relations involving conjugation by matrices in
$GL(n; {\mathbb C})$ and direct sums. Quite generally, on such a matricially generalized open set $\O$, the
correspponding matricially generalized scalar analytic functions form a noncommutative algebra $A(\O)$ and
there is a generalization $\p$ of the difference quotient derivation-comultiplication which yields a topological
GDQ ring structure. In the $C^*$-context, if $\O=\O^*$ in a suitably defined sense, $A(\O)$ becomes a
$*$-algebra and there is also a notion of dual positivity.

The duality map appears as a transformation from $E'$ to an $A(\O)$, where $\O$ is the matricially generalized
resolvent set and this transformation intertwines GDQ ring structures and positivity on $E'$ with
dual-positivity on $A(\O)$.

Besides section 1, which is the introduction, there are eight more sections numbered 2 to 9.

Section 2 contains preliminaries on GDQ rings.

Section 3 is about the new corepresentations we found.

In section 4 we introduce multivariable GDQ rings and we give a reduction result to a one-variable GDQ ring in
case $n=p^2, n$ the number of "variables``.
We also prove a result about how $\p_{X:B}$ arises in general GDQ rings.

Section 5 deals with full $B$-resolvents and resolvent sets, which are the matricial $B$-valued generalizations
of usual resolvents and resolvent sets.

Section 6 takes up the matricial generalization of functions and sets which go with the generalized resolvents.

Section 7 gives the construction of the topological GDQ ring structure on the algebras $A(\O)$ of fully
matricial functions. We have prefered to define the derivation-comultiplication as taking values in some
"two-variable`` $A(\O;\O)$ instead of entering here the technical problems about tensor-products and topologies
on the $A(\O)$'s.

Section 8 contains a discussion of dual positivity in $A(\O)$.

Section 9 introduces the duality $U$-transform and discusses its intertwining properties for GDQ structure and
positivity.

{\bf Acknowledgments}: Most of this work was done while the author held an International Blaise Pascal Research
Chair, from the State and Ile de France Region, managed by the Fondation de l'Ecole Normale Sup\'{e}rieure,
 and visited the Institut de Math\'{e}matiques de Jussieu during the Spring 2003. He was also supported in part
by NSF Grant DMS-0079945.

\section{Preliminaries on GDQ Rings}
\label{sec2}

We shall slightly amend the terminology in [12] by breaking up the definition of a GDQ ring with 
involution into smaller groups of conditions.\\

{\bf 2.1. Definition.}
{\em A generalized difference quotient ring (a $GDQ$ ring) is an object $(A, \mu, \p)$, where $A$ is an algebra
over
${\mathbb C}$ and}\\
\\
(GDQ 1) \ \ \ \ \ $\mu: A\otimes A\rightarrow A$ {\em is the multiplication map}\\
\\
(GDQ 2) \ \ \ \ \ $\p:A\rightarrow A\otimes A$ {\em is linear and coassociative, i.e.} $(\p \ox {\rm id}_A)
\circ \p =  ({\rm id}_A \ox \p) \circ \p$\\
\\
(GDQ 3) \ \ \ \ \ $\p$ {\em is a derivation, i.e.} $\p \circ \mu = ({\rm id}_A \ox \mu) \circ (\p \ox {\rm id}_A) +
(\mu \ox {\rm id }_A) \circ ({\rm id}_A \ox \p) .$\\
\\
{\em In general we do not require A to have a unit. If $1\in A$ is a unit then the GDQ ring will be called unital.
}\\

{\bf 2.2. Remark.} {\em A GDQ ring can always be made unital by adjoining a unit $1$ and putting $\p 1 = 0$.}
\\

{\bf 2.3. Definition.} {\em $(A, \mu , \p , L)$ is a graded GDQ ring if $(A, \mu , \p)$ is a GDQ ring 
and there is a linear map $L: A \rightarrow A $ (the grading), so that\\
\\
{\rm (L1)} \ \ \ \ \ $L - {\rm id}_A$ is a derivation of the algebra $(A, \mu)$\\
\\
{\rm (L2)} \ \ \ \ \ $L$ is a coderivation of the coalgebra $(A, \p)$, i.e. $ \p \circ  L = (L \ox {\rm id}_A + 
{\rm id}_A \ox L)\circ \p$}.
\\

{\bf 2.4 Definition} {\em An involution of a GDQ ring $(A, \mu, \p)$ is a conjugate-linear involution 
$A \ni a \rightarrow a^* \in A$ of the vector space $A$ so that \\
\\
{\rm (I1)} \ \ \ \ \ $(A, \mu , *)$ is an algebra with involution\\
\\
{\rm (I2)} \ \ \ \ \ $\p(a^*) = \s_{12}((\p a)^*)$, where $*$ on $A\ox A$ is given by $(x\ox y)^* = x^*\ox y^*$ 
and $\s_{12}(x\ox y) = y\ox x$.\\

If $L$ is a grading, compatibility with the involution means}\\
\\
(I3) \ \ \ \ \ $L(a^*) = (L(a))^*$.\\

If $V$ is a vector space, we shall denote by $V^\bullet$ its dual endowed with the topology of pointwise convergence.
By $\widehat{\ox}$ we denote the projective tensor product. The duality theorem (Thm. 5.3 in [12]) can be
restated  in the following form.\\

{\bf 2.5. Theorem.} {\em If $(A, \mu . \p)$ is a GDQ ring then $(A^\bullet, \p^\bullet , \mu^\bullet )$ 
satisfies the GDQ ring conditions with $\ox$ replaced by $\widehat{\ox}$. If $L$ is a grading and $*$ is an
involution for $(A, \mu . \p)$, then $\L = - {\rm id}_{A^\bullet} + L^\bullet $ and $\xi^*(a) =
\overline{\xi(a^*)}$ satisfy the grading and respectively the involution conditions for $(A^\bullet ,
\p^\bullet ,
\mu^\bullet)$.} (with $\ox$ replaced by $\widehat{\ox}$).\\

By ${\mathfrak M}_p(A)$ we denote the $p\x p$ matrices over $A$, an individual matrix being either written in
the form $(a_{ij})_{1\leq i, j \leq p}$ or $\sum_{1 \leq i, j \leq p} a_{ij} \ox e_{ij}$ where $a_{ij} \in A$ 
and $e_{ij}$ are the matrix-units. A corepresentation of $(A, \mu , \p)$ is a matrix $\a = \sum_{i, j}a_{ij}
\ox e_{ij} \in {\mathfrak M}_p(A)$ so that 
\[
 \sum_{i,k} \p a_{ik} \ox e_{ik} = \sum_{i,j,k} a_{ij} \ox a_{jk} \ox e_{ik}.
\]
This can also be written: 
\[
(\p \ox {\rm id}_{{\mathfrak M}_p})\a = \a  \ox_{{\mathfrak M}_p} \a.
\]
The main result about corepresentations (Prop. 1.4 in [12]) is the following.\\

{\bf 2.6. Theorem. } {\em Let $(A, \mu , \p)$ be a unital GDQ ring and assume $X\in A$ is so  that
 $\p X = 1 \ox 1$. If $\a \in {\mathfrak M}_p(A)$ is invertible, the following are equivalent:

{\rm (i)}  $\a$ is a corepresentation

{\rm (ii)} $\a = ((n_{ij} - X\d_{ij})^{-1})$ where $n_{ij} \in N = {\rm Ker} \d$. }\\

Since this is a functional analysis paper, the algebraic facts will guide our functional analysis constructions,
even if they are not directly applicable. This is a familiar situation from the theory of Kac algebras and 
$C^*$-quantum groups, where finding the appropriate topological tensor products and topological duals are 
subtle analysis questions.

In particular the vague idea, that the dual object should be constructed by mapping $\var \in A^\bullet$ into 
the direct sum of $(\var (a_{ij}))_{1\leq i, j \leq p} \in {\mathfrak M}_p$ where
 $\a = (a_{ij})_{1 \leq i,j \leq p}$ runs over a sufficiently large set of corepresentations of $(A, \mu , \p)$
 poses many analytical problems.

\section{More Corepresentations}
\label{sec3}
Throughout this section $(A, \mu , \p)$ will denote a unital GDQ ring and $X \in A$ will be an element so that
 $\p X = 1\ox 1$. We shall exhibit corepresentations which enlarge the set provided by Theorem. 2.6.

Like in [12], it will be convenient to use $d:{\mathfrak M}_p(A) \rightarrow {\mathfrak M}_p(A\ox A)$ and 
$d = \p \ox {\rm id}_{{\mathfrak M}_p}$ which is a derivation w.r.t. the bimodule structure given by the
homomorphisms $\var_1, \var_2 : {\mathfrak M}_p(A) \rightarrow {\mathfrak M}_p(A \ox A)$ so that 
$\var_1(a_{ij})_{1 \leq i, j \leq p} = (a_{ij}\ox 1)_{1\leq i,j \leq p}$, 
$\var_2(a_{ij})_{1 \leq i, j \leq p} = (1\ox a_{ij})_{1\leq i,j \leq p}$. We shall also denote $X\ox I_p$
by $x\in {\mathfrak M}_p(A)$ and write $1$ for the unit $1\ox I_p$ of ${\mathfrak M}_p(A)$.\\

{\bf 3.1. Proposition.} {\em If $N \in {\rm Ker} J$ and $\b_1, \b_2, \b_3 \in {\mathfrak M}_p(N)$ are such that
 $\b_2 - \b_1(X \ox I_p)\b_3$ is invertible, then
\[
\a = \b_3(\b_2 - \b_1(X\ox I_p)\b_3)^{-1}\b_1
\]
is a corepresentation.\\

Proof.} Let $y = \b_2 - \b_1(X\ox I_p)\b_3 = \b_2 - \b_1 x \b_3$. We have
\[
d(y^{-1}) = -\var_1(y^{-1})d(y)\var_2(y^{-1})
\]
Hence
\begin{eqnarray*}
d(\a) & = & d(\b_3y^{-1}\b_1) = \var_1(\b_3)d(y^{-1})\var_2(\b_1) \\
      & = & -\var_1(\b_3)\var_1(y^{-1})(-\var_1(\b_1)\var_2(\b_3))\var_2(y^{-1})\var_2(\b_1)\\ 
      & = & \var_1(\a)\var_2(\a)
\end{eqnarray*}

 which is the desired result.
\qed

We have also the following general procedure for producing more corepresentations.

{\bf 3.2. Lemma.} {\em Let $\xi \in {\mathfrak M}_p(A)$ be a corepresentation and let $\b \in {\mathfrak M}_p({\rm
Ker} \p)$. If $1- \xi\b$ is invertible then $\a =(1 - \xi\b)^{-1}\xi$ is a corepresentation.

If $1-\b\xi$ is invertible then $\g = \xi(1-\b\xi)^{-1}$ is a corepresentation.\\

Proof.} Because of symmetry, we will only prove the first assertion.\\
We have
\begin{eqnarray*}
d((1 - \xi\b)^{-1}\xi) &=& -\var_1((1 - \xi\b)^{-1})d((1 - \xi\b))\var_2((1 - \xi\b)^{-1})\var_2(\xi) + 
      \var_1((1 - \xi\b)^{-1})\var_1(\xi)\var_2(\xi)\\
&=&\var_1((1 - \xi\b)^{-1})\var_1(\xi)\var_2(\xi)\var_2(\b)\var_2((1 - \xi\b)^{-1})\var_2(\xi) +
       \var_1((1 - \xi\b)^{-1})\var_1(\xi)\var_2(\xi)\\
&=& \var_1((1 - \xi\b)^{-1}\xi)\var_2(\xi + \xi \b (1 - \xi\b)^{-1}\xi)\\
&=& \var_1(\a)\var_2(\xi + ((1 - \xi\b)^{-1} -1)\xi )\\
&=& \var_1(\a)\var_2(\a).
\end{eqnarray*}
\qed

\section{Reduction of Multivariable GDQ Rings}
\label{sec4}

Studying $\p_{X:B}$ does not mean a limitation to one variable. In this section we briefly explain how 
multivariable situations can easily be reduced to the $\p_{X:B}$ setting.\\

{\bf 4.1.} The typical multivariable situation deals with $A = B\<X_1, \dots, X_n\>$, the ring of 
noncommutative polynomials in the noncommutative variables $X_1,\dots ,X_n$ and with noncommutative
scalars $B$.
This means monomials are of the form $b_0X_{i_1}b_1X_{i_2}b_2\dots X_{i_n}b_n$ and the only relations are
 those arising from $1X_j = X_j1 = X_j$. The $n$ partial difference quotients $\p_i :A\rightarrow A\ox A$ are 
the derivations such that $\p_i X_j = \d_{ij}1\ox 1$ and $\p_iB = 0$. Thus each $(A, \mu ,\p_i)$ is a GDQ ring 
and we have compatibility relations
\[
(\p_i\ox {\rm id}_A) \circ \p_j = ({\rm id}_A \ox \p_j) \circ \p_i
\]
for all $1 \leq i, j \leq n$. In particular
\[
(A, \mu , \l_1\p_1 + \dots + \l_n\p_n)
\]
is again a GDQ ring ($\l_j \in {\mathbb C}$). Note that for a ``multivariable GDQ ring" 
$(A, \mu, \p_1, \dots , \p_n)$, i.e. the structure in which $(A,\mu , \p_j)$ are GDQ rings and the 
compatibility relations hold, $(A, \mu. \sum_i\l_{i1}\p_i, \dots , \sum_i\l_{in}\p_i)$ is again a
multivariable GDQ.\\

{\bf 4.2} In case $n = p^2$ a multivariable GDQ ring $(A, \mu, \p_1,\dots , \p_{p^2})$ can be replaced
by a one-variable $(\tilde{A},\tilde{\mu}, \p)$. More precissely, we take $\tilde{A} = {\mathfrak M}_p(A) =
{\mathfrak M}_p\ox A$ (where ${\mathfrak M}_p$ is short for ${\mathfrak M}_p({\mathbb C})$). We may reindex $\p_1, 
\dots , \p_{p^2}$ (possibly preceded by a linear transformation, if we want to preserve some involution) and replace
 them by $\p_{ij}, 1\leq i,j\leq p$. Let further
\[
\D_{ij} = \sum_{1\leq k\leq p}e_{ki}\ox e_{jk} \in {\mathfrak M}_p \ox {\mathfrak M}_p
\]
Note that:
\[
(T \ox I_p)\D_{ij} = \D_{ij}(I_p \ox T)
\]
if $T \in {\mathfrak M}_p$.

We then define
\[
\p :\tilde{A} \rightarrow \tilde{A} \ox \tilde{A}
\]
by
\[
\p(T\ox a) =
 \sum_{1\leq i,j\leq p}((T \ox I_p)\D_{ij})\ox \p_{ij}a \in({\mathfrak M}_p)^{\ox 2} \ox A^{\ox 2}\simeq \tilde{A}^{\ox 2}
\]
where the isomorphism takes $(T_1\ox T_2)\ox (a_1 \ox a_2)$ to $(T_1\ox a_1)\ox(T_2\ox a_2)$.

Note that we also have
\[
\p(T\ox a) = \sum_{1\leq i, j\leq p}(\D_{ij}(I_p\ox T))\ox \p_{ij}a.
\]

That $\p$ is a derivation is seen by the computation
\begin{eqnarray*}
\p((T_1\ox a_1)(T_2\ox a_2)) &=& \p(T_1T_2\ox a_1a_2) = \sum_{i,j}((T_1T_2\ox I_p)\D_{ij})\ox \p_{ij}a_1a_2\\
&=&\sum_{i,j}((T_1\ox I_p)\D_{ij}(I_p\ox T_2))\ox((a_1\ox 1)(\p_{ij}a_2) + (\p_{ij}a_1)(1\ox a_2))\\
&=& ((T_1\ox a_1) \ox (I_p\ox 1))\p(T_2\ox a_2) + (\p(T_1 \ox a_1))((I_p \ox 1)\ox(T_2 \ox a_2)).
\end{eqnarray*}

Before checking coassociativity remark that
\[
\p(e_{rs}\ox a) = \sum_{i, j}(e_{ri}\ox e_{js}) \ox(\p_{ij}a)
\]

We have:
\begin{eqnarray*}
((\p \ox {\rm id})\circ \p)(e_{rs}\ox a) &=& \sum_{i,j}(\p \ox {\rm id})((e_{ri}\ox e_{js}) \ox \p_{ij}a)\\
&=&  \sum_{k,l} \sum_{i,j} (e_{rk}\ox e_{li}\ox e_{js}) \ox((\p_{kl} \ox {\rm id})\circ \p_{ij})a)
\end{eqnarray*}
while on the other hand
\begin{eqnarray*}
(({\rm id} \ox \p) \circ \p)(e_{rs}\ox a) &=& \sum_{k.l}({\rm id} \ox \p)((e_{rk} \ox e_{ls}) \ox \p_{kl}a)\\
&=& \sum_{i,j}\sum_{k,l}(e_{rk}\ox e_{li} \ox \e_{js})\ox((({\rm id} \ox \p_{ij})\circ \p_{kl})a)
\end{eqnarray*}
and the coassociativity follows from compatibility of $\p_{ij}$ and $\p_{kl}$.
\\

{\bf 4.3} If there are elements $Y_{ij}\in A$ so that $\p_{rs}Y_{ij} = \d_{ri}\d_{sj}1\ox 1$ then it is 
easily seen that
\[
Y = \sum_{1\leq i, j \leq p}e_{ij} \ox Y_{ij} \in \tilde{A}
\]
will have the property $\p Y = (I_p\ox 1)\ox(I_p \ox 1)$. Also if
\[
Z = \sum_{1\leq i,j\leq p}e_{ij}\ox Z_{ij}
\]
then $\p Z = 0$ is equivalent to $\p_{rs}Z_{ij} = 0$ for all $1\leq r, s,i.j\leq p$,
so that
\[
{\rm Ker} \p = {\mathfrak M}_p \ox (\cap_{1\leq i,j\leq p} {\rm Ker}\p_{ij}).
\]\\

{\bf 4.4} Returning to the multivariable GDQ ring $A=B\<X_1, \dots , X_{p^2}\>$ and the partial free difference 
quotients $\p_1,\dots , \p_{p^2}$ w.r.t. $X_1,\dots , X_{p^2}$, the preceding construction, combined with a 
linear transformation, gives the following. We consider $\tilde{A} = 
{\mathfrak M}_p \ox B\<X_1,\dots , X_{p^2}\>$ which is isomorphic to $D\<X\>$ where $D= {\mathfrak M}_p \ox B$
and $X = \sum_jT_j \ox X_j$ for some basis $T_j, 1\leq j\leq p^2$ of ${\mathfrak M}_p$. The replacement for the 
multivariable GDQ ring is then $D\<X\>$ with comultiplication-derivation $\p_{X:D}$.

Note that in case $B = {\mathbb C}$ or $B = {\mathfrak M}_q$, this reduction has the pleasant feature that
$D$,  which is ${\mathfrak M}_p$ or ${\mathfrak M}_{pq}$, is finite-dimensional.\\

{\bf 4.5.} We conclude section 4 with a structure result for GDQ rings, which was used implicitely in the
preceding  subsection.\\

{\bf Proposition.} {\em Let $(A,\mu,\p)$ be a GDQ ring with unit and assume there is $X\in A$ such that 
$\p X = 1\ox 1$. Let further $N={\rm Ker}\p$.

Then the canonical homomorphism $\psi:N\<X\> \rightarrow A$ is an injection and it is a GDQ ring homomorphism
when $N\<X\>$ is endowed with the comultiplication $\p_{X:N}$.\\

Proof.} A derivation being completely determined by the way it acts on the generators of an algebra, the only
 assertion we really need to prove is the injectivity of $\psi$. Let \mbox{$\psi_k:N^{\ox (k+1)}\rightarrow A$}
be the linear maps so that
\[
\psi_k(n_0\ox \dots \ox n_k) = n_0Xn_1X\dots n_k.
\] 
We must prove ${\rm Ker} \psi_k = 0$ and the ranges of the $\psi_k, (k\geq 0)$ are linearly independent.
Iterating $\p$ we define $\p^{(k)}= (\p\ox {\rm id}_{k-1}) \circ \p^{(k-1)},\ \p^{(1)} = \p$. Then
$\p^{(k)}\psi_k (N^{\ox (k+1)}) \subset N^{\ox (k+1)}$ and $N^{\ox(k+1)}\;|\; \p^{(k)}\circ \psi_k = {\rm
id}_{N^{\ox(k+1)}}$ and $ \p^{(k)}\circ \psi_l =0$ if $l < k$. The assertion follows from these facts.
\qed

\section{The full $B$-Resolvent}
\label{sec5}

Let $E$ be a Banach algebra with unit, let $1\in B\subset E$ be a closed subspace containing the unit and let
 $Y \in E$ be an element. The concepts we examine in this section will serve also as motivating examples
 in the next and it is good to note that the case when $B$ is a Banach subalgebra is of particular interest.\\

{\bf 5.1. Definition.} {\em The set
\[
\rho_n(Y;B) = \{b\in {\mathfrak M}_n(B)\ |\ Y\ox I_n - b \ \ {\rm invertible}\}
\]
will be called the n-th B-resolvent set of $Y$. The collection of subsets $(\rho_n(Y;B))_{n\geq1}$ of
the ${\mathfrak M}_n(B)$'s will be called the full B-resolvent of $Y$. The operator-valued function 
$R_N(Y;B)(\cdot): \rho_n(Y;B) \rightarrow {\mathfrak M}_n(E)$ defined by $R_n(Y;B)(b) = (Y\ox I_n -
b)^{-1}$
 will be called the n-th B-resolvent of Y. The collection of functions $(R_n(Y;B))_{n\geq 1}$ will be called
 the full B-resolvent of $Y$.}\\

Some basic facts about these concepts are summarized in the next proposition.\\

{\bf 5.2. Proposition.} 
\\
(i)   $\rho_n(Y;B)$ {\em is open in} ${\mathfrak M}_n(B)$.
\\
\\
(ii) $\rho_m(Y;B) \oplus \rho_n(Y;B) = \rho_{m+n}(Y;B) \cap ({\mathfrak M}_m(B)\oplus {\mathfrak M}_n(B))$.
\\
\\
(iii) $(S\ox 1)\rho_n(Y;B)(S\ox 1)^{-1} = \rho_n(Y;B)$ {\em if} $S\in GL(n; {\mathbb C})$.
\\
\\
(iv) {\em If $b' \in \rho_m(Y;B)$, $b'' \in \rho_n(Y;B)$ and if $\b \in {\mathfrak M}_{m,n}(B)$ is an $m \x n$
 matrix with entries in B, then}
\[
\left(\begin{array}{cc}
b'& \b\\
0 & b''
\end{array}\right) \in \rho_{m+n}(Y;B).
\]
\\
\\
(v)   $R_n(Y;B)$ {\em is a complex analytic function.}
\\
\\
(vi)  {\em If $b' \in \rho_m(Y;B)$, $b''\in \rho_n(Y;B)$ then $R_{m+n}(Y;B)(b'\oplus b'') =
 R_m(Y;B)(b')\oplus R_n(Y;B)(b'')$.}
\\
\\
(vii) {\em If $b\in {\mathfrak M}_n(B)$ and $S\in GL(n; {\mathbb C})$ then
\[
R_n(Y;B)((S\ox 1)b(S\ox 1)^{-1})= (S\ox 1)R_n(Y;B)(S\ox 1)^{-1}.\]
 
Proof.} Most assertions are rather obvious and will be left to the reader. We will only prove (iv). In view of
 (ii), $b\oplus b' \in \rho_{m+n}(B)$ and in view of (i)
\[
\left(\begin{array}{cc}
b' & \e\b\\
0 & b''
\end{array}\right) \in\rho_{m+n}(B)
\]
for some $\e \neq 0$. Then (iv)  follows by applying (iii) with $S = \e^{- 1}I_m\oplus I_n$.
\qed

\section{Fully Matricial Functions and Sets}
\label{sec6}
$G$ and $H$ will be Banach spaces over ${\mathbb C}$. If $S\in GL(n; {\mathbb C})$ and $T\in {\mathfrak M}_n$
 we denote by $Ad S$ the automorphism of  ${\mathfrak M}_n$ so that $(AdS)(T) = STS^{-1}$. The corresponding
automorphism of ${\mathfrak M}_n(H) = {\mathfrak M}_n\ox H$ will be denoted by $AdS\ox I_{\cal H}$ or simply
$AdS\ox I$ and its action is $(AdS\ox I)(T\ox h) = STS^{-1}\ox h$.\\

{\bf 6.1. Definition.} {\em A fully matricial G-set is a sequence $(\O_n)_{n\geq 1}$ so that}
\\

(FMS1)\ \ \ \ \ $\O_n \in {\mathfrak M}_n(G)$
\\

(FMS2)\ \ \ \ \ $\O_{m+n} \cap ({\mathfrak M}_m\oplus {\mathfrak M}_n)= \O_m\oplus \O_n$
\\

(FMS3)\ \ \ \ \ $(AdS\ox I)(\O_n) = \O_n$ {\em if} $S\in GL(n; {\mathbb C})$.
\\

{\em A fully matricial G-set is open or closed if each $\O_n$ is open or respectively closed.}\\

{\bf 6.2. Proposition.} {\em If $(\O_n)_{n\geq 1}$ is a fully matricial open set and if $g'\in \O_m,\;g'' \in
\O_n$ and $\g \in {\mathfrak M}_{m,n}(G)$, then} $\left(\begin{array}{cc} g' & \g\\0 & g''\end{array}\right)
\in \O_{m+n}$.\\

The proof is along the same lines as the proof of (iv) in Proposition 5.2.

In case $G = {\mathbb C}$, using the Jordan form of a matrix it is possible to describe the fully matricial 
${\mathbb C}$-sets.\\

{\bf 6.3. Proposition.} (i) {\em A fully matricial ${\mathbb C}$-set $(\O_n)_{n\geq 1}$ is described in an unique
 way by giving for each $\l \in {\mathbb C}$ an additive subsemigroup $L(\l )\subset {\mathbb N}$. Then
$T\in \O_n$ iff for each eigenvalue $\l \in \s(T)$, the length of the corresponding Jordan blocks in the 
Jordan form of $T$ are in $L(\l )$.}

(ii) {\em $(\O_n)_{n\geq 1}$ is a closed (respectively open) fully matricial ${\mathbb C}$-set iff $\O_1$ is 
closed (respectively open) and $\O_n = \{T \in {\mathfrak M}_n\ |\ \s(T) \subset \O_1 \}$. (In particular if 
the fully matricial ${\mathbb C}$-set is closed or open, the $L(\l)$'s can only be $\emptyset$ or ${\mathbb N}$}.
\\

The proof of (i) is an exercise in combining the Jordan form with the similarity and direct sum properties 
of fully matricial sets, which we leave to the reader. We will only explain the different reasons in (ii) when 
$\O_n$ is closed or open, why the $L(\l)$'s can only be ${\mathbb N}$ or $\emptyset$. In both cases, using
(FMS3) the disscussion breaks down to showing that if $T \in \O_n$ is an upper triangular matrix, then its 
(1,1)-entry $\l_1$ will be in $\O_1$.

If the fully matricial set is closed let $S(\l)$ be the diagonal matrix with entries $1,\l,\dots , \l$.
Then
\[
\lim_{\l\rightarrow \i}S(\l)TS(\l)^{-1} = T'
\]
where $T'$ is the direct sum of the $1\x 1$ matrix $\l_1$ and an $(n-1)\x (n-1)$ matrix. Since $\O_n$ is closed
 $T'\in \O_n$ and by (FMS2), $\l_1 \in \O_1$.

If $\O_n$ is open, we can find $T' \in \O_n$ so that 
$T' = \left(\begin{array}{cc}\l_1 & *\\0 & S\end{array}\right)$ where $S \in {\mathfrak M}_{n-1}$ is so that
 $\s(S) \not\ni\l_1$. Then using (FMS3) and the Jordan form we find $\left(\begin{array}{cc}\l_1 & 0\\0 &
S\end{array}\right) 
\in \O_n$ and by (FMS2) $\l_1 \in \O_1$.\\

{\bf 6.4. Proposition.} {\em If $(\O_n^{(i)})_{n \geq 1}$, $(i\in I)$ are fully matricial G-sets, then 
$(\bigcap_{i \in I}\O_n^{(i)})_{n\geq 1}$is a fully matricial G-set.}\\

The proof is left as an exercise.

{\em In particular the family of open fully matricial G-sets is stable under such finite componentwise intersections.

Similarly, the family of closed fully matricial G-sets is stable under arbitrary componentwise intersections.}

It seems natural to consider the topology (viewed for instance as subsets of $\coprod_{n\geq 1}{\mathfrak
M}_n(G)$) generated by the open fully matricial $G$-sets.\\

{\bf 6.5. Definition.} {\em A fully matricial H-valued function on a fully matricial G-set $(\O_n)$ is a 
sequence $(R_n)_{n\geq 1}$ so that}
\\

(FMF1)\ \ \ \ \ {\em $R_n:\O_n \rightarrow {\mathfrak M}_n(H)$ is a function}
\\

(FMF2)\ \ \ \ \ {\em If $g'\in \O_m\:, g'' \in \O_n$ then $R_{m+n}(g'\oplus g'') = R_m(g') \oplus R_n(g'')$}
\\

(FMF3)\ \ \ \ \  {\em If $S\in GL(n; {\mathbb C})$ and $g\in \O_n$ then $R_n((AdS\ox I_G)(g)) = 
(AdS\ox I_H)(R_n(g))$.} 
\\

{\em A fully matricial function is continuous if each component is continuous. A fully matricial function 
is analytic if the fully matricial G-set on which it is defined is open and the components $R_n$ are analytic.}\\

{\bf 6.6. Remark.} A fully matricial function amounts to a sequence of functions the graphs of 
which form a fully matricial $G\x H$-set.\\

{\bf 6.7. Lemma.} {\em Let $(R_n)_{n\geq 1}$ be a continuous fully matricial H-valued function on the fully 
matricial G-set $(\O_n)_{n\geq 1}$. Assume $g' \in \O_m\:, g''\in \O_n$ and $\g \in {\mathfrak M}_{m, n}(G)$.
Then for some $h \in {\mathfrak M}_{m,n}(H)$
\[
 R_{m+n}(\left(\begin{array}{cc}g' & \l \g\\0 & g''\end{array}\right)) = \left(\begin{array}{cc}R_m(g') & \l h\\0 & R_n(g'')\end{array}\right)
\]
for all $\l \in {\mathbb C}$.\\

Proof.} Let
$
R_{m+n}(\left(\begin{array}{cc}g' & \g\\0 & g''\end{array}\right)) = 
\left(\begin{array}{cc}h' & h_{12}\\h_{21} & h''\end{array}\right)$ and let $S(\e) = \e I_m\oplus I_n \in 
GL(m+n; {\mathbb C})$. Then if $\l \not= 0$
\begin{eqnarray*}
R_{m+n}(\left(\begin{array}{cc}g' & \e \l \g\\0 & g''\end{array}\right)) & = & R_{m+n}((AdS(\e \l)\ox I_G)
\left(\begin{array}{cc}g' & \g\\0 & g''\end{array}\right))\\
&=&(AdS(\e \l)\ox I_H)\left(\begin{array}{cc}h' & h_{12}\\h_{21} & h''\end{array}\right)\\
   &=&\left(\begin{array}{cc}h' & \e \l h_{12}\\\l^{-1} \e^{-1}h_{21} & h''\end{array}\right)\,.
\end{eqnarray*}
\\
Since $R_{m+n}$ is continous and $lim_{\e \rightarrow 0}\left(\begin{array}{cc}g' & \e \g\\0 & g''\end{array}\right)
= \left(\begin{array}{cc}g' & 0\\0 & g''\end{array}\right) \in \O_{m+n}$ we infer
\[
lim_{\e \rightarrow 0}\left(\begin{array}{cc} h' & \e \l h_{12}\\\l^{-1}\e^{-1}h_{21} & h''\end{array}\right) = 
\left(\begin{array}{cc}R_m(g' )& 0\\0 & R_n(g'')\end{array}\right)
\]
and hence $h_{21} = 0, h' = R_m(g'), h''= R_n(g'')$.
\qed

{\bf 6.8. Remark.} The reader has probably recognized by now that the full B-resolvent set and the
 full $B$-resolvent are examples of a fully matricial $B$-set and respectively of an analytic $B$-valued fully 
matricial function defined on an open fully matricial $B$-set.\\

{\bf 6.9. Definition.} {\em A fully matricial G-set $(\O_n)_{n\geq 1}$ will be called finite if it also satisfies
\[
\left(\begin{array}{cc}g' & \g\\0 & g''\end{array}\right) \in \O_{m+n} \Rightarrow g'\in \O_m\,, g''\in \O_n\,.
\]}
\\

{\bf 6.10. Remark.} Full resolvent sets in a finite von Neumann algebra provide examples of finite 
fully matricial sets. On the other hand, taking $E= B= {\cal B}(l^2({\mathbb N})), Y=0, m=n=1$,
the full $B$-resolvent of $Y$ is not a finite fully matricial $B$-set since 
$\left(\begin{array}{cc}S & K\\0 & S^*\end{array}\right)$, where $S$ is the unilateral and $K$ a 
rank one operator making the matrix to be the bilateral shift, is in $\rho_{1+1}(0;B)$ without $S, S^*$ being 
in $\rho_1(0;B)$. 
\\

{\bf 6.11} Returning to the context of Proposition 6.4, we associate with an open fully matricial $G$-set
$\O = (\O_n)_{n\geq 1}$ the set $\tilde{\O} = \coprod_{n\geq 1}\O_n \subset \coprod_{n\geq 1} {\mathfrak M}_n(G)$
and we consider the topology they generate on $\coprod_{n\geq 1}{\mathfrak M}_n(G)$. It is then also natural to
associate with
 $\tilde{\O}$ the analytic or continuous fully matricial $H$-valued functions on $\O$ and the sheaves on
 $\coprod_{n\geq 1}{\mathfrak M}_n(G)$ which they generate.
\\

{\bf 6.12. Abbreviations.}  From now on we will also use the abbreviations {\bf FM G-S} for fully matricial 
$G$-set and {\bf FMF} for fully matricial $G$-function. Also  FMS will abreviate fully matricial set
 and FMAF will abbreviate fully matricial analytic function.

\section{The GDQ Ring of Scalar Fully Matricial Analytic Functions}
\label{sec7}

{\bf 7.1.} Let $\O =(\O_n)_{n\geq 1}$ be an open FM G-S. To avoid amending our assumptions on $G$ to introduce
 more structure, we will assume that $G$ {\em is an operator system} (see [3]) i.e. it is isomorphic to a space 
of operators on Hilbert space which is selfadjoint and unital and is correspondingly endowed with involution, 
unit and is matrix-normed. (The reader could simplify and assume $G$ is a unital $C^*$-algebra.)  We should
also
 clarify from the beginning that the term GDQ ring in the title of this section has been used rather loosely:
the tensor product required for the comultiplication would be a topological one, and we would actually circumvent
this question interpreting the tensor product as some two-variable function. Our aim here is to clarify the
 function-theory aspect of the comultiplication and to return to precise topological  GDQ ring structure later.\\

{\bf 7.2.} Let $A(\O)$ denote the ${\mathbb C}$-valued FMAF on $\O$. If $r=(r_n)_{n\geq 1}, s=(s_n){n\geq 1}$
are  in $A(\O)$, then $r+s=(r_n+s_n)_{n\geq 1}$ and $rs=(r_ns_n)_{n\geq 1}$ are in $A(\O)$ which thus is
naturally a
 noncommutative ring. Moreover $1=(1_n)_{n\geq 1}$, where $1_n$ denotes the constant function on $\O_n$ with 
value the identity $n\x n$ matrix, is the unit in $A(\O)$.

Let $\O^* = (\O_n^*)_{n\geq 1}$, where $\O_n^* = \{T^* \:| \: T\in \O_n\}$. If $r\in A(\O )$ we define
 $r^* = (r_n^* )_{n\geq 1} \in A(\O)$ by $r_n^*(g) = (r_n(g^*))^*$. Thus $r\rightarrow r^*$ is a conjugate-linear 
antiisomorphism of $A(\O)$ and $A(\O^*)$. {\em In case $\O = \O^*$ this makes $A(\O)$ a unital algebra 
with involution.}\\

{\bf 7.3.} Let $K=(K_n)_{n\geq 1}$ be a fully matricial subset of $\O$. We will say {\em $K$ is properly
included in $\O$ if
\[
\sup_{n\in {\mathbb N}} \sup_{k \in K_n}\| k\|_n < \i
\]
and if there is $\e > 0$ such that $K_n + \e({\mathfrak M}_n(G))_1 \subset \O_n$ for all $n\in {\mathbb N}$.}
(Here $\|\ \|_n$  is the norm and $({\mathfrak M}_n(G))_1$ the unit ball in ${\mathfrak M}_n(G)$.
 Clearly this definition uses the fact that $G$ is matrix-normed). If $r\in A(\O)$ we define
\[
\|r\|_K = \sup_{n\in {\mathbb N}} \sup_{k_n\in K_n}\|r_n(k_n)\|_n
\]
where $\|\cdot \|_n$ is the norm on ${\mathfrak M}_n({\mathbb C})$. {\em Unless $\|r\|_K < \i$, for all
properly  included $K$, it may be natural to add this condition and consider the corresponding subalgebra
$A_{pr}(\O)$ of
$A(\O)$.}\\

{\bf 7.4} The comultiplication derivation will be defined piecewise, i.e. for fixed matrix-sizes. We 
will use algebras of matrix-valued analytic functions $A_{n_1,\dots ,n_p}(\O_{n_1}; \dots ; \O_{n_p})$, where 
$n = n_1+ \dots + n_p$, consisting of analytic maps $f\rightarrow \O_{n_1} \x\dots \x \O_{n_p} \rightarrow 
{\mathfrak M}_{n_1} \ox \dots \ox {\mathfrak M}_{n_p}$ which are $GL(n_1) \x \dots \x GL(n_p)$-equivariant
\[
f((AdS_1 \ox I_G)(g^{(1)}), \dots , (AdS_p \ox I_G)(g^{(p)})) =
((AdS_1)\ox \dots \ox (AdS_p))f(g^{(1)},\dots , g^{(p)})
\]
where $S_j \in GL(n_j), g^{(j)} \in \O_{n_j}$.

A result similar to Lemma 6.7 holds for functions in $A_{m+n}(\O_{m+n})$.\\

{\bf 7.5. Lemma.} {\em Let $f_{m +n} \in A_{m+n}(\O_{m+n})$. Then:}\\

(i) {\em if $g'\in \O_m, g''\in \O_n$ there are $a'
\in {\mathfrak M}_m, a'' \in {\mathfrak M}_n$ so that
\[ f_{m+n}(g' \oplus g'') = \ a'\oplus a''\].}\\

(ii) {\em if $g', g'', a', a''$ are as in (i) and $\p \in {\mathfrak M}_{m,n}(G)$ there is 
$h\in {\mathfrak M}_{m,n}$ so that \[
f_{m+n}(\left(\begin{array}{cc}g' & \l \g\\0 & g''\end{array}\right)) =
 \left(\begin{array}{cc}a' & \l h\\0 & a''\end{array} \right)
\]
for all $\l \in {\mathbb C}$.\\

Proof. } (i) If $f_{m+n}(g'\oplus g'') = \left( \begin{array}{cc}h' & h_{12}\\h_{21} & h''\end{array}\right )$ 
then in view of the equivariance applied to $S=\e I_m \oplus I_n$ we get $h_{12} = \e h_{12}, h_{21} = 
\e^{-1}h_{21}$, so that $h_{12} = 0, h_{21} = 0$.\\

(ii) If
\[
f_{m+n}(\left( \begin{array}{cc}g' & \g\\0 & g''\end{array}\right )) = \left( \begin{array}{cc}h' & h_{12}\\h_{21} & h''\end{array}\right )
\]
conjugation with $\e I_M \oplus I_n$ yields
\[f_{m+n}(\left( \begin{array}{cc}g' & \e \g\\0 & g''\end{array}\right ))=
\left( \begin{array}{cc}h' & \e h_{12}\\\e^{-1} 
h_{21} & h''\end{array}\right )
\] and since
\[
\lim_{\e\rightarrow 0}\left( \begin{array}{cc}h' & \e h_{12}\\\e^{-1} h_{21} & h''\end{array}\right ) =
\left(\begin{array}{cc}a' & 0\\0 & a''\end{array}\right)
\]
we infer $h'=a', h''=a'', h_{21} =0$. Hence
\[
f_{m+n}(\left( \begin{array}{cc}g' & \e \g\\0 & g''\end{array}\right)) = \left(\begin{array}{cc}a' &
 \e h_{12}\\0 & a'' \end {array}\right).
\]
\qed

{\bf 7.6} There is a canonical identification $\a$ of ${\mathfrak M}_m \ox {\mathfrak M}_n$  with the linear
 operators ${\cal L}({\mathfrak M}_{mn})$ on ${\mathfrak M}_{mn}$. If $a\in {\mathfrak M}_m, b\in {\mathfrak M}_n, 
c\in {\mathfrak M}_{mn}$
\[ (\a(a\ox b))(c) = acb \in {\mathfrak M}_{mn}.
\]
If $m, n$ need to be specified we will write $\a_{mn}$.\\

{\bf 7.7.} We define 
\[
\p_{m,n} : A_{m+n}(\O_{m+n})\rightarrow  A_{m,n}(\O_m ;\O_n)
\]
as follows. Let $f\in A_{m+n}(\O_{m+n}), h\in {\mathfrak M}_{m,n}, 1\in G, g'\in \O_m, g''\in \O_n$ and let
 $\g_{m,n} :{\mathfrak M}_{m,n} \rightarrow {\mathfrak M}_{m+n}$ be the map which puts ${\mathfrak M}_{m,n}$
 into the right $m\x n$ corner of ${\mathfrak M}_{m+n}$ i.e. $\g_{m,n}(e_{jk}) = e_{j, m+k}$ and $\g_{mn}$
 is linear. By lemma 7.5 
\[
\frac{{\rm d}}{{\rm d}\e}f(g'\oplus g'' +\e \g_{m,n}(h)\ox 1)|_{\e = 0} = \g_{m,n}(h')
\]
for some $h'\in {\mathfrak M}_{m,n}$. Hence for each $(g', g'')\in \O_m \ox \O_n$ we get a map
 ${\mathfrak M}_{m,n} \ni h \rightarrow h' \in {\mathfrak M}_{m,n}$. Applying $\a^{-1}$ to this map gives an 
element in ${\mathfrak M}_m\ox {\mathfrak M}_n$ which is our definition of $(\p_{m,n}f)(g', g'')\in
 {\mathfrak M}_m \ox {\mathfrak M}_n$. This can also be written as a formula. Since the differential of $f$
 at $g'\oplus g''$ is a linear map, we have
\[
\p_{m,n}f(g', g'') = \sum_{\begin{array}{c}1\leq i, j\leq m\\1\leq k, l\leq n \end{array}}
(\frac{{\rm d}}{{\rm d}\e}f(g'\oplus g'' +\e e_{j,k+m}\ox 1)|_{\e = 0})_{i, m+l} e_{ij}^{(m)}\ox e_{kl}^{(n)}
\]
where $(\cdot)_{i, m+l}$ denotes the $(i, m+l)$ entry of the $(m+n)\x(m+n)$ matrix.

It is clear that $\p_{m,n}f$ defined in this way is an analytic function $\O_m \x \O_n \rightarrow 
{\mathfrak M}_m \ox {\mathfrak M}_n$.\\

{\bf 7.8.} To check that $\p_{m,n}f$ is a $GL(m)\x GL(n)$ equivariant map, remark first that if $S'
\in GL(m), S'' \in GL(n)$ then
\begin{eqnarray*}
\lefteqn{\frac{\rm d}{{\rm d}\e}f((AdS'\ox I_G)g'\oplus (AdS''\ox I_G)g'' + 
\e\g_{m, n}(S'hS''{}^{-1}) \ox 1)|_{\e = 0} =} \\
& = & Ad(S'\oplus S'')(\frac{\rm d}{{\rm d} \e}f(g'\oplus g'' +\e\g_{m,n}(h)\ox 1))|_{\e =0} =\\
&= & S'\g_{m, n}(h')S''{}^{-1}.
\end{eqnarray*}
Thus we must check that, if $T\in {\cal L}({\mathfrak M}_{m,n})$ and $\tilde{T} \in {\cal L}({\mathfrak M}_{m,n})$ 
is given by $\tilde{T}(S'hS''{}^{-1}) = S'(T(h))S''{}^{-1}$, then $\a^{-1}(\tilde{T}) =
 ((AdS')\ox (AdS''))\a^{-1}(T)$. This is the same as the following equivariance for $\a$: if $T = \a(\xi)$ and 
$\tilde{T} = \a((AdS' \ox AdS'')\xi)$ then $\tilde{T}(S'hS''{}^{-1}) = S'T(h)S''{}^{-1}$. It suffices to see this 
for $\xi = a\ox c$. Then $T(h) = ahc$ and $\tilde{T}(h) = (S'aS''{}^{-1})h(S''cS''{}^{-1})$ so that
$\tilde{T}(S'hS''{}^{-1}) = S'ahcS''{}^{-1} = S'T(h)S''{}^{-1}$, i.e. the equivariance we wanted to check. Hence
$\p_{m,n}f_{m+n}\in A_{m,n}(\O_m;\O_n)$.

The derivation property will be obtained from the following lemma.\\

{\bf 7.9. Lemma.} {\em If $f, \tilde{f} \in A_{m+n}(\O_{m+n})$ and $g'\in \O_m, g'' \in \O_n, a', b' \in 
{\mathfrak M}_m, a'', b'' \in {\mathfrak M}_n, f(g'\oplus g'') = a' \oplus a'', \tilde{f}(g'\oplus g'') =
 b'\oplus b''$ then 
\[
(\p_{m,n}(f\tilde{f}))(g', g'') = (a' \ox I_n)(\p_{m,n}\tilde{f})(g', g'') + (\p_{m,n}f)(g',g'')(I_m\ox b'').
\]

Proof.} In view of Lemma 7.5 we have
\[
f(g'\oplus g'' + \l \g_{m,n}(h)\ox 1) = f(g'\oplus g'') + \l \frac{{\rm d}}{{\rm d}\e} f(g'\oplus g'' + 
\e \g_{m,n}(h) \ox 1)|_{\e = 0}
\]
and the same holds with $f$ replaced by $\tilde{f}$. Multiplying we get
\begin{eqnarray*}
(f\tilde{f})(g'\oplus g'' + \l \g_{m,n}(h)\ox 1) &=& 
 f(g'\oplus g'')\tilde{f}(g'\oplus g'') + \\ &+&\l f(g'\oplus g'')\frac{{\rm d}}{{\rm d} \e}\tilde{f}
(g'\oplus g'' + \e \g_{m,n}(h)\ox 1)|_{\e = 0}+\\&+& \l \frac{{\rm d}}{{\rm d} \e}f(g'\oplus g''
+ \e \g_{m,n}(h)\ox 1)|_{\e = 0}\tilde{f}(g'\oplus g'') \\ &+& O(\l^2).
\end{eqnarray*}
This gives
\begin{eqnarray*}
\frac{\rm d}{{\rm d}\e}(f\tilde{f})(g'\oplus g'' + \e \g_{m,n}(h)\ox 1)|_{\e = 0} &=& (a'\oplus a'')
\frac{{\rm d}}{{\rm d} \e}\tilde{f}(g'\oplus g'' + \e \g_{m,n}(h)\ox 1)|_{\e = 0} + \\
&+& \frac{\rm d}{{\rm d} \e}f(g'\oplus g'' + \e \g_{m,n}(h)\ox 1)|_{\e = 0}(b' \oplus b'').
\end{eqnarray*}
Taking the right $(m,n)$-block corner gives
\[
\a(\p_{m,n}(f\tilde{f})(g'g''))(h) = a'\a(\p_{m,n}(\tilde{f})(g',g''))(h) + \a(\p_{m,n}(f)(g',g''))(h)b''
\]
The result follows from\\

$a'\a(\xi)(h) = \a((a'\ox I_n)\xi)(h)$\\

$\a(\xi)b''h = \a(\xi(I_m\ox b''))(h)$\\
\\
if $\xi \in {\mathfrak M}_m\ox {\mathfrak M}_n$.
 \qed

{\bf 7.10 Corrolary.} {\em If $r=(r_n)_{n\geq 1} \in A(\O), s=(s_n)_{n\geq 1} \in A(\O )$ and $g'\in \O_m, g''\in 
\O_n$ then
\begin{eqnarray*}
(\p_{m,n}(rs)_{m+n})(g', g'') &=& (r_m(g')\ox I_n)(\p_{m,n}s_{m+n})(g', g'')\\&+&
(\p_{m,n}r_{m+n})(g',g'')(I_m 
\ox s_n(g'')).
\end{eqnarray*}
}
 This is immediate from the preceding Lemma when we take into account that $r_{m+n}(g'\oplus g'')=
 r_m(g')\oplus r_n(g'')$ and $s_{m+n}(g'\oplus g'')= s_m(g')\oplus s_n(g'')$.\\

{\bf 7.11.} To combine the $\p_{m,n}$ into a derivation for $A(\O)$ we will need to define "several variables
fully matricial analytic functions".

Let $\O^{(j)},
 j=1,\dots , p$ be FM G-S. 

{\em We define the p-variables scalar fully matricial analytic functions on $\O^{(1)}\x \dots \x \O^{(p)}$ to
 be families of analytic functions $(f_{n_1, \dots , {n_p}})_{{n_1}\geq 1, \dots , {n_p}\geq 1}$
where $f_{n_1, \dots , {n_p}}: \O^{(1)}_{n_1} \x \dots \ox \O^{(p)}_{n_p} \rightarrow {\mathfrak M}_{n_1} \ox
\dots
\ox {\mathfrak M}_{n_p}$ which are $GL(n_1) \x \dots \x GL(n_p)$-equivariant and so that

\begin{eqnarray*}
f_{n_1,\dots ,n_{j-1},n'_j +n''_j,n_{j+1},\dots, n_p}(g_1,\dots ,g'_j\oplus g''_j,\dots g_p) &=&
f_{n_1,\dots ,n_{j-1},n'_j,n_{j+1},\dots, n_p}(g_1,\dots , g'_j,\dots , g_p)\\
&\oplus & f_{n_1,\dots ,n_{j-1},n''_j,n_{j+1},\dots, n_p}(g_1,\dots ,g''_j,\dots , g_p).
\end{eqnarray*}
The scalar p-variables FMAFs on $\O^{(1)}\x \dots \x\O^{(p)}$ will be denoted by $A(\O^{(1)}; \dots ;\O^{(p)})$.
Clearly, $A(\O^{(1)}; \dots ;\O^{(p)})$ is an algebra with unit.

If $f\in A(\O^{(1)}; \dots ;\O^{(p)})$ and $\tilde{f}\in A(\tilde{\O}^{(1)}; \dots ;\tilde{\O}^{(q)})$ then we
define \\
$f\ox \tilde{f} \in A(\O^{(1)}; \dots ;\O^{(p)};\tilde{\O}^{(1)}; \dots ;\tilde{\O}^{(q)})$ by\\
\[(f\ox \tilde{f})_{n_1,\dots , n_p, \tilde{n}_1,\dots ,\tilde{n}_q}(g_1,\dots ,g_p, \tilde{g}_1,\dots ,
\tilde{g}_q) = f_{n_1,\dots ,n_p}(g_1,\dots ,g_p)\ox \tilde{f}_{\tilde{n}_1,\dots ,\tilde{n}_q}(\tilde{g}_1, 
\dots ,\tilde{g}_q).\]}

{\bf 7.12. Lemma.} {\em If $r\in A(\O)$ then $(\p_{m,n}r_{m+n})_{m\geq 1, n\geq 1} \in A(\O;\O)$.

Proof.} Analyticity and equivariance have already been checked and we are left with
$\p_{m,n}r_{m+n}(g'\oplus g'', g)=
(\p_{{m_1},n}r_{{m_1}+n})(g',g)\oplus (\p_{{m_2},n}r_{{m_2}+n})(g'',g)$ where $g'\in \O_{m_!}, g''\in \O_{m_2},
g\in~\O_n$, $m=m_1+m_2$ and $(\p_{m,n_1+n_2}r_{m+n})$ where now $g\in \O_m, g'\in \O_{n_1}, g''\in \O_{n_2}, n=
n_1+n_2$. The direct sums are in the sense of
\[
 ({\mathfrak M}_{m_1}\ox {\mathfrak M}_n)\oplus ({\mathfrak M}_{m_2} \ox {\mathfrak M}_n)=({\mathfrak M}_{m_1}
\oplus {\mathfrak M}_{m_2})\ox {\mathfrak M}_n \subset {\mathfrak M}_{m_1+m_2} \ox{\mathfrak M}_n\]
and in the second case of
\[
({\mathfrak M}_m\ox {\mathfrak M}_{n_1})\oplus ({\mathfrak M}_m\ox {\mathfrak M}_{n_2}) = {\mathfrak M}_m
\ox ({\mathfrak M}_{n_1}\oplus {\mathfrak M}_{n_2}) \subset {\mathfrak M}_m\ox {\mathfrak M}_{n_1+n_2}.\]
We will only sketch how one checks the first of the two equalities for $\p_{m,n}$, the second one being similar.

First remark that $\a^{-1}$ behaves well w.r.t direct sums, i.e. if $T_1\in{\cal L}({\mathfrak M}_{{m_1}, n}),
T_2
\in {\cal L}({\mathfrak M}_{{m_2}, n})$ and $T_1\oplus T_2 \in{\cal L}({\mathfrak M}_{{m_1},n}\oplus
 {\mathfrak M}_{{m_2}, n}) = {\cal L}({\mathfrak M}_{{m_1}+{m_2}, n})$ then
$\a^{-1}_{{m_1}+{m_2},n}(T_1\oplus T_2) = \a^{-1}_{{m_1},n}(T_1)\oplus \a^{-1}_{{m_2},n}(T_2)$.
Thus it will suffice to check that
\[
(\a_{{m_1}+{m_2},n}\d_{m,n}r_{m,n})(g'\oplus g'',g) = (\a_{{m_1},n}\p_{{m_1},n}r_{{m_1},n})(g',g) \oplus
(\a_{{m_2},n}\p_{{m_2},n}r_{{m_2},n})(g',g)\]
In view of Lemma 7.5 and of the direct sum property of a FMAF, what we must prove amounts to the following. 
Let $h_1,h_1' \in {\mathfrak M}_{{m_1},n}$ and $h_2, h_2' \in {\mathfrak M}_{{m_2},n}$ be such that
\begin{eqnarray*}
r_{{m_1},n}(\left(\begin{array}{cc}g' & h_1\ox 1 \\0 & g\end{array}\right)) &=& 
\left(\begin{array}{cc}a' & h_1' \\0 & a\end{array}\right)\\
r_{{m_2},n}(\left(\begin{array}{cc}g'' & h_2\ox 1 \\0 & g\end{array}\right)) &=&
\left(\begin{array}{cc}a'' & h_2' \\ 0 & a\end{array}\right)
\end{eqnarray*}
(where $a' = r_{m_1}(g'), a''= r_{m_2}(g''), a = r_n(g)$). Then we will have
\[
r_{{m_1}+{m_2}+n }(\left(\begin{array}{ccc} g'& 0 & h_1\ox1\\0  & g''  &h_2\ox 1 \\0 & 0 & g\end{array}\right))
= \left(\begin{array}{ccc} a' & 0  & h_1' \\ 0 & a''  &h_2' \\ 0& 0 & a\end{array}\right) .
\]
Since
\[
r_{{m_1}+{m_2}+n}(\left(\begin{array}{ccc} g'& 0 & h_1\ox1\\0  & g''  &h_2\ox 1 \\0 & 0 & g\end{array}\right))
=\left(\begin{array}{cc} r_{{m_1}+{m_2}}(\left(\begin{array}{cc} g' & 0  \\ 0 & g'' \end{array}\right))
 &  * \\ 0 & r_n(g) \end{array}\right)
\]
by Lemma 7.5 and $r_{{m_1}+m_2}(g'\oplus g'') = a'\oplus a''$, all we need to check is that the (1,3) and
(2,3) block-entries of the result are $h_1'$ and $h_2' $. This can be done by several applications of 
direct sum and $GL$-equivariance properties:
\begin{eqnarray*}
r_{m_1+n+m_2+n}(
\left(\begin{array}{cccc} g' & h_1\ox1  &0  &0 \\ 0 &  g &0 &0 \\0 &0  &g'' &h_2\ox 1 \\0 &0 &0 & g
\end{array}\right)) &=&
\left(\begin{array}{cccc} a' & h_1'  & 0 & 0\\ 0 & a  &0 & 0\\0 & 0 &a'' &h'_2 \\0 &0 &0 &a\end{array}\right)\\
\\
\\
r_{m_1+m_2+n+n}(
\left(\begin{array}{cccc} g' &  0 &h_1\ox 1  &0 \\ 0 & g''  &0 &h_2\ox 1 \\0 & 0 &g &0 \\0 &0 & 0&g
\end{array}\right))&=&
\left(\begin{array}{cccc} a' & 0  & h_1' &0 \\ 0 & a''  &0 &h_2' \\0 & 0 &a &0 \\0 &0 &0 &a\end{array}\right)\\
\\
\\
r_{m_1+m_2+n+n}(
\left(\begin{array}{cccc} g' & 0  &h_1\ox 1  &0 \\0  & g''  &h_2\ox 1 &h_2\ox 1 \\ 0& 0 &g &0 \\0 &0 &0 &g
\end{array}\right)&=&
\left(\begin{array}{cccc} a' & 0  & h_1' &0 \\ 0 &a'' & h_2'&h_2' \\0 & 0 &a &0 \\ 0&0 &0 &a\end{array}\right)\\
\end{eqnarray*}
on one hand and also on the other hand
\[
=\left(\begin{array}{cc} r_{m_1+m_2+n}(\left(\begin{array}{ccc}  g'&  0 & h_1\ox 1 \\0  &g''&h_2\ox 1 \\
0 & 0 &g \end{array}\right))& *\\ 0 & a \end{array}\right)
\]
which then gives the desired result.
\qed

{\bf 7.13.} {\em If $r\in A(\O)$ we shall denote by $\p r$ the element $(\p_{m,n}r_{m+n})_{m\geq 1, n\geq 1}
\in A(\O ;\O)$.} 

Before going further let us also record the following fact which appeared in the preceding proofs.

{\bf Lemma.} {\em Let $g'\in \O_m, g''\in \O_n,h,h'\in {\mathfrak M}_{m,n}$ be such that
\[
r_{m+n}(\left(\begin{array}{cc} g' & h\ox 1  \\ 0 & g'' \end{array}\right)) =
\left(\begin{array}{cc} r_m(g') &  h' \\ 0 &r_n(g'') \end{array}\right)
\]
Then $h' = (\a_{m,n}(\p_{m,n}r_{m+n})(g', g''))(h)$.

In particular the map taking $h$ to $h'$ is linear and takes $sht$ to $sh't$ if $s\in GL(m), t\in GL(n)$.}\\

{\bf 7.14.} We pass to the coassociativity property of $\p$. Since we have not identified $A(\O; \O)$ with
 a tensor product $A(\O)\ox A(\O)$ we will define maps $({\rm id}\ox \p)A(\O; \O) \rightarrow A(\O;\O ;\O)$
and respectively $(\p \ox {\rm id})A(\O; \O) \rightarrow A(\O;\O ;\O)$. The most convenient seems to be to use
the
 formula for matrix entries given at the end of section 7.7. Thus, {\em we define for $h\in
A_{m,n+p}(\O;\O_{n+p})$
 and $g\in \O_m, g'\in \O_n, g''\in \O_p$
\begin{eqnarray*}
\lefteqn{
(({\rm id} \ox \p)_{m,n,p}k)(g, g', g'')=}\\
& & \sum_{\begin{array}{c}1\leq a, b\leq m\\1\leq c, d\leq n\\1\leq e,f\leq p \end{array}}(\frac{\rm d}{{\rm
d}
\e}k(g; g'\oplus g'' +\e e_{d, n+e})|_ {\e = 0})_{(a,b)(c, n+f)}e^{(m)}_{ab}\ox e^{(n)}_{cd}\ox e^{(r)}_{ef}
\end{eqnarray*}
where the index $(ab)(c, n+f)$ stands for the coefficient of $e^{(m)}_{ab}\ox e^{(n+p)}_{c, n+f}$.

 In particular if $f\in A_m(\O_m), \tilde{f}\in A_{n+p}(\O_{n+p})$ then $({\rm id} \ox \p)_{m,n,p}
(f\ox \tilde{f}) = f\ox (\p_{m,n}\tilde{f})$.

We leave it to the reader to check that $({\rm id}\ox \p)_{m,n,p}k \in A_{m+n+p}(\O_m; \O_n; \O_p)$. Part of
the verification can be done using $g\in \O_m, g'\in \O_n, g''\in \O_p$, functionals $\var \in
 ({\mathfrak M}_m)'$ , the functions $(\var \ox {\rm id})k(g;\cdot)\in A(\O_{n+p})$ , the fact that
$
((\var\ox {\rm id}_{{\mathfrak M}_n}\ox {\rm id}_{{\mathfrak M}_p})({\rm id}\ox \p)_{m,n,p}k)(g; g', g'')=
\p_{n,p}((\var \ox {\rm id}_{{\mathfrak M}_{n+p}})k(g; \cdot))(g', g'')
$
and the results we already have for $\p_{n,p}$. Using this type of argument one then checks that if $k =
 (k_{n_1, n_2})_{n_1\geq 1, n_2\geq 1}\in A(\O;\O)$ then $(({\rm id}\ox \s)_{n_1, n_2, n_3}
k_{n_1,n_2+n_3})_{n_1\geq 1, n_2\geq 1, n_3\geq 1} \in A(\O; \O ;\O)$ and a similar result for $\p \ox {\rm
id}$.

Checking that $({\rm id}\ox \p)\circ \p = (\p \ox {\rm id}) \circ \p$, after we've pushed aside all these
 questions boils down to the following result.}\\

{\bf 7.15 Lemma} {\em If $k\in A_{m+n+p}(\O_{m+n+p})$ then $({\rm id}\ox  \p)_{m,n,p}\p_{m,n+p}k =
(\p \ox {\rm id})_{m,n,p}\p_{m+n,p}k$.\\

Proof.} Let $g\in \O_m, g'\in \O_n, g''\in \O_p$. We have
\begin{eqnarray*}
\lefteqn{
(({\rm id}\ox \p)_{m,n,p} \circ \p_{m,n+p}k)(g,g',g'')_{(a,b)(c,d)(e,f)} =}\\
&=& \frac{{\rm d}}{{\rm d}\e_2}(
\frac{{\rm d}}{{\rm d}\e_1}(k(g\oplus g'\oplus g''+ \e_1e_{b,m+c}+\e_2e_{m+d,m+n+e}))_{a,m+n+f}|_{\e_1=0})|_
{\e_2 = 0}
\end{eqnarray*}
On the other hand
\begin{eqnarray*}
\lefteqn{
((\p\ox{\rm id})_{m,n,p} \circ \p_{m+n,p}k)(g,g',g'')_{(a,b)(c,d)(e,f)} =}\\
&=& \frac{{\rm d}}{{\rm d}\e_2}(
\frac{{\rm d}}{{\rm d}\e_1}(k(g\oplus g'\oplus g''+ \e_1e_{m+d,m+n+e}+\e_2e_{b,m+c}))_{a,m+n+f}|_{\e_1=0})|_
{\e_2 = 0}
\end{eqnarray*}
The equality of the two quantities is thus quite obvious.\\

{\bf 7.16. Lemma.} {\em Let $\O_1\subset {\mathbb C}$ be an open set, $G= {\mathbb C}$ and $\O_n =
 \{a\in {\mathfrak M}_n \; |\; \s(a)\subset \O_1\}$ . Let further $f=(f_n)_{n\geq 1}\in A(\O)$ where 
$\O =(\O_n)_{n\geq 1}$, so that $f_n(a) = f_1(a)$, where the right hand side has the meaning of functional
 calculus. Then if $z_1, z_2 \in \O_1, z_1 \not= z_2,$
\[
(\p_{1,1}f_2)(z_1, z_2) = \frac{f_1(z_1)-f_1(z_2)}{z_1-z_2}.
\]

Proof.} Let $(\p_{1,1}f_2)(z_1,z_2) = \l\in {\mathfrak M}_1^{\ox 2} \simeq{\mathbb C}$. Then
\[
\left(\begin{array}{cc} f_1(z_1) &\l \\ 0 & f_1(z_2) \end{array}\right)=
f_1(\left(\begin{array}{cc} z_1 &  1 \\0 & z_2 \end{array}\right)).
\]
Since $\left(\begin{array}{cc} z_1 & 1 \\ 0 &z_2 \end{array}\right) = (Ad
\left(\begin{array}{cc}1  &  (z_2-z_1)^{-1} \\ 0 & 1 \end{array}\right))
\left(\begin{array}{cc} z_1 & 0  \\ 0 &z_2  \end{array}\right)$ it follows that
\[
\left(\begin{array}{cc} f_1(z_1) &\l \\ 0 & f_1(z_2) \end{array}\right)=
Ad \left(\begin{array}{cc} 1 &  (z_2-z_1)^{-1} \\ 0 & 1 \end{array}\right)\left(\begin{array}{cc} f_1(z_1) & 1
\\ 0 &f_2(z_2) \end{array}\right)
\]
which gives the desired result.
\qed

\section{Dual Positivity in $A(\O)$}

Let $\O$ be an open FM G-S over the operator space $G$ and assume $\O=\O^*$. We will use the map $\a\p =\nabla$.
In particular if $g'\in \O_m, g''\in \O_n$ then $(\n_{m,n}f_{m+n})(g', g'')$ is an element in ${\cal
L}({\mathfrak M}_{m,n})$.\\

{\bf 8.1. Definition.} {\em $f\in A(\O)$ is dual positive if $f = f^*$ and for any $g\in \O_n, n\in {\mathbb
N}$,
\[(\n_{m,n}f)(g,g^*):{\mathfrak M}_n\rightarrow {\mathfrak M}_n
\]
is a positive map (i.e. transforms positive operators into positive operators).}\\

{\bf 8.2. Proposition.} {\em If $f\in A(\O)$, the following are equivalent:}\\

(i) {\em $f$ is dual positive}\\

(ii) {\em $f=f^*$ and for any $g^{(j)}\in \O_{n(j)}, 1\leq j \leq p, \oplus_{1\leq i,j \leq p}(\n_{n(i), n(j)}f)
(g^{(i)}, g^{(j)*})$ is a positive linear map of $\oplus_{i,j}{\mathfrak M}_{n(i),n(j)}$, identified with 
${\mathfrak M}_{n(1)+\dots +n(p)}$, into itself}\\

(iii) {\em $f=f^*$ and for any $g\in \O_n$ the map $(\n_{n,n}f)(g,g^*): {\mathfrak M}_n \rightarrow 
{\mathfrak M}_n$ is completely positive.\\

Proof.} Clearly (ii) $\Rightarrow$ (i) and (iii)$\Rightarrow$(i).

 (i)$\Rightarrow $(ii). It suffices to show 
that the map $\oplus_{1\leq i,j\leq p}(\n_{n(i),n(j)}f_{n(i)+n(j)})(g^{(i)}, g^{(j)*})$ coincides with the map
\[
(\n_{n,n}f_{n+n})(g^{(1)}\oplus \dots \oplus g^{(p)}, g^{(1)*}\oplus \dots \oplus g^{(p)*})
\]
where $n = n(1) + \dots + n(p)$. Indeed in view of the definition of $\a$, this is the same as establishing that
$\oplus_{1\leq i,j \leq p}\p_{n(i),n(j)}f_{n(i)+n(j)}(g^{(i)}, g^{(j)*})$ as an element of 
$\oplus_{1\leq i,j\leq p}{\mathfrak M}_{n(i)}\ox {\mathfrak M}_{n(j)}~\subset~{\mathfrak M}_n \ox {\mathfrak
M}_n$ coincides with $\p_{n,n}f_{n+n}(g^{(1)}\oplus \dots \oplus g^{(p)},g^{(1)*}\oplus \dots \oplus g^{(p)*})
\in {\mathfrak M}_n\ox {\mathfrak M}_n$. This in turn is an immediate consequence of the fact that $\p f \in
A(\O;\O)$.

(ii) $\Rightarrow$ (iii) If $t^{(ij)} \in {\mathfrak M}_n, 1\leq i,j\leq p$ form a $p\x p$ matrix with $n\x n$
blocks, which is positive in ${\mathfrak M}_{np}$, we must show that the $np\x np$ matrix formed from the blocks 
$(\n_{n,n}f_{n+n}(g,g^*))(t^{(ij)})$ is also positive. this is precisely the statement in (ii) in case 
$n(1) = \dots = n(p) = n$ and $g^{(1)} = \dots g^{(p)} = g$.
\qed

\section{The Full Resolvent Transform U}
\label{sec9}

{\bf 9.1} The dual GDQ ring corresponds to a map of the dual of the GDQ ring into a GDQ ring of the $A(\O)$
type. As long as we don't use an involution we will stay in the context of section 5. Thus, $E$ will be a Banach
algebra with unit,
$1\in B\subset E$ a Banach subalgebra and $Y\in E$ an element. Let $\rho(Y;B) = (\rho_n(Y;B))_{n\geq 1}$ be the
the full $B$-resolvent set of $Y$ and $R(Y;B)= (R_n(Y;B))_{n\geq 1}$  the full $B$-resolvent.

{\em By ${\cal RA}(Y;B)$ we shall denote the subalgebra of $E$ generated by $B, \{Y\}$ and the matrix coefficients 
of the $\{R_n(Y;B)(b)\;|\; n\in {\mathbb N}, b\in \rho_n(Y;B)\}$.}\\

{\bf 9.2.} {\em We shall assume there is a derivation-comultiplication
\[
\p: {\cal RA}(Y;B) \rightarrow {\cal RA}(Y;B) \ox {\cal RA}(Y;B)
\]
so that ${\cal RA}(Y;B)$ is a GDQ ring and $\p B = 0, \p Y = 1\ox 1$. If such a $\p$ exists, then it is unique,}
i.e. it is completely determined by the conditions $\p B = 0, \p Y = 1\ox 1$. Indeed, the $R_n(Y;B)(b)$ will 
then be corepresentations, and the corresponding equation determines $\p$ on the matrix coefficients. Thus
 $\p$ is completely determined on the generators of ${\cal RA}(Y;B)$, hence being a derivation it is
completely determined on ${\cal RA}(Y;B)$.\\

{\bf 9.3.} {\em We shall also assume ${\cal RA}(Y;B)$ is dense in $E$.} 

Let ${\cal CR}(Y;B)$ denote the matrix coefficients of $R_n(Y;B)(b)$ ($b\in \rho_n(Y;B), n\in {\mathbb N}$).\\

{\bf 9.4. Lemma.} {\em ${\cal CR}(Y;B)$ is closed under multiplication. The assumptions in 9.3 imply that the
linear span of ${\cal CR}(Y;B)$ is dense in $E$.

Proof.} Remark first that if $a\in {\mathfrak M}_m(E), \tilde{a} \in {\mathfrak M}_n(E), x\in 
{\mathfrak M}_{m,n}(E)$ and $a^{-1}, \tilde{a}^{-1}$ exist, then
\[
\left(\begin{array}{cc} a & -x  \\ 0 & \tilde{a} \end{array}\right)^{-1} = 
\left(\begin{array}{cc} a^{-1} &a^{-1}x\tilde{a}^{-1} \\ 0 &\tilde{a}^{-1}  \end{array}\right)
\]
In particular the $(i, l+m)$-entry of this $2\x 2$ block matrix is the $(i,l)$-entry of $a^{-1}x\tilde{a}^{-1}$.
Choosing $x$ to be the $(j,k)$ matrix unit we find that for this choice of $x$ one of the matrix coefficients
of $\left(\begin{array}{cc} a & -x  \\ 0 & \tilde{a} \end{array}\right)^{-1}$ is the product of the $(i,j)$-entry 
of $a^{-1}$ and of the $(k,l)$-entry of $\tilde{a}^{-1}$. Taking $a$ and $\tilde{a}$ to be $\b-Y\ox I_m$ 
and respectively $\b'-Y\ox I_n$ we get that ${\cal CR}(Y;B)$ is closed under multiplication.

 Thus the linear span of ${\cal CR}(Y;B)$ is an algebra and to prove the second assertion it suffices to prove 
that its closure contains $Y$ and $B$. Since the linear span of invertible elements
in $B$ is $B$ it will suffice 
to prove that the invertible elements in $B$ and $Y$ are in the closure of the linear span of ${\cal CR}(Y;B)$.
If $b\in B$ is invertible, then so is $(\l b-Y)$ for $\l$ large enough and $(\l b-Y)^{-1}\in{\cal CR}(Y;B)$ and
$\lim_{\l
\rightarrow
\i}\l^{-1} (\l b- Y)^{-1} = b^{-1}$. The assertion about $Y$ follows from
\[
Y=\lim_{\e\rightarrow 0}\e^{-1}(\e^{-1}(\e^{-1}-Y)^{-1} -\e^{-2}(\e^{-2} - Y)^{-1})
\]
\qed

{\bf 9.5.} Let $E'$ denote the dual of the Banach space $E$. The {\em full resolvent transform is defined to be
the map
\[
U:E'\rightarrow A(\rho((Y;B))
\]
so that $U(\var)=(U_n(\var))_{n\geq 1}$ where 
\[
U_n(\var)(\cdot) =(\var\ox {\rm id}_{{\mathfrak M}_n})(R_n(Y;B)(\cdot)) \in A_n(\rho_n(Y;B)).
\]
}
(Remark that $U_n(\var)(\cdot) $ is fully matricial analytic because $R_n(Y;B)(\cdot)$ is fully matricial analytic.)
\\

{\bf 9.6. Proposition.} {\em If $\var_1, \var_2, \var_3 \in E'$ are such that 
\[
\var_1(a) = (\var_2 \ox \var_3)(\p a)
\]
for all $a\in {\cal RA}(Y;B)$, then 
\[U(\var_1) = U(\var_2)U(\var_3)
\]\\

Proof.} It is actually sufficient that the assumption hold for $a\in {\cal CR}(Y;B)$ in order to get the 
conclusion. Indeed applying the assumption to each matrix coefficient of $R_n(Y;B)(b) =\a$ we have that 
\begin{eqnarray*}
U_n(\var_1)(b) &=& (\var_1 \ox {\rm id}_{{\mathfrak M}_n})(\a) = (\var_2 \ox \var_3 \ox {\rm id}_{{\mathfrak M}_n}
)(\p\a)\\
&=& (\var_2 \ox \var_3 \ox {\rm id}_{{\mathfrak M}_n}
)(\a \ox_{{\mathfrak M}_n} \a)\\&=&
U_n(\var_2)(\a)U_n(\var_3)(\a) 
\end{eqnarray*}
\qed

{\bf 9.7} Before stating the duality property involving the comultiplication of $A(\rho (Y;B))$ we need to clarify 
a notation we'll use. If $b= \sum_{i,j}b_{i,j}\ox e_{i,j}^{(m)} \in {\mathfrak M}_m(B)$ and $b' = \sum_{k,l}b'_{k,l}
\ox e_{kl}^{(n)}\in {\mathfrak M}_n(B)$ we denote by $b\ox_B b' \in {\mathfrak M}_{mn}(B)$ the $mn \x mn$ matrix,
 or equivalently the element in ${\mathfrak M}_m\ox {\mathfrak M}_n\ox B$ given by 
$\sum_{i,j,k,l}e_{ij}^{(m)}\ox e_{kl}^{(n)} \ox(b_{ij}b'_{kl})$. Equivalently if $a\ox \b \in {\mathfrak M}_m\ox
B$ and $a'\ox \b' \in {\mathfrak M}_n\ox B$ then $(a\ox \b)\ox_B(a'\ox \b') = a\ox a'\ox \b\b'$.\\

{\bf 9.8. Proposition.} {\em If $\var \in E'$ and $b_1 \in \rho_m(Y;B), b_2\in \rho_n(Y;B)$ then
\[
(\var \ox {\rm id}_{{\mathfrak M}_m}\ox {\rm id}_{{\mathfrak M}_n})(R_m(Y;B)(b_1)\ox_E R_n(Y;B)(b_2))=
 -\p_{m,n}(U_{m+n}(\var ))(b_1; b_2).
\]

Proof.} Returning to the computations on which the proof of Lemma 9.4 relies, let $a= b_1 - Y\ox I_m,
 \tilde{a}= b_2 - Y\ox I_n$ and let $x=1\ox e_{j,k}^{(m,n)}$ so that
\[
(R_m(Y;B)(b_1)\ox_E R_n(Y;B)(b_2))_{(i,j)(k,l)} = (\left(\begin{array}{cc} a &-x   \\ 0 &\tilde{a}
  \end{array}\right)^{-1})_{i,l+m}.
\]
On the other hand 7.7 and Lemma7.13 combined give that
\begin{eqnarray*}
(\p_{m,n}(U_{m+n}(\var))(b_1;b_2))_{(i,j)(k,l)}&=& ( -U_{m+n}(\var )
\left(\begin{array}{cc} b_1 & -x  \\ 0 & b_2\end{array}\right))_{i,m+l}\\
&= & -((\var\ox {\rm id}_{{\mathfrak M}_{m+n}})\left(\begin{array}{cc} a & -x  \\0  &\tilde{a} 
\end{array}\right)
 ^{-1})_{i,m+l}
\end{eqnarray*}
which implies the desired result.\qed

{\bf 9.9.Remark.} Propositions 9.6 and 9.8 express the fact that the $U$-transform relates "dual GDQ structure"
on 
$E'$ with the "topological GDQ structure" of $A(\rho(Y;B))$ endowed with the comultiplication $-\p$. These
duality statements take this indirect form because of the rather algebraic setting of our discussion
 (i.e. without analytic assumptions on the comultiplication of ${\cal RA}(Y;B)$ and a closer examination of the
 topological tensor product in the GDQ structure of $A(\rho(Y;B))$ ).\\

{\bf 9.10.Proposition.} 

(i) {\em U is injective.}\\

(ii) {\em $\var \in E'$ satisfies the trace-condition $\var([E,E])=0$ iff $\p_{m,n}(U_{m+n}(\var))(b_1; b_2)=\\
= \e\circ\p_{n,m}(U_{m+n}(\var))(b_2;b_1)$ for all $b_1\in \rho_m(Y;B), b_2\in \rho_n(Y;B),m\geq 1, n\geq 1$.}\\
(here $\e:{\mathfrak M}_m\ox {\mathfrak M}_n \rightarrow {\mathfrak M}_n\ox {\mathfrak M}_m$ permutes the two
factors)\\

 {\em Proof.} (i) follows from Lemma 9.4.

(ii) If $\var ([E,E])= 0$ then the equality we want to prove, in view of Proposition 9.8 is equivalent to
\[
\var([(R_m(Y;B)(b_1))_{i,j},(R_n(Y;B)(b_2))_{k,l}]) = 0
\]
which follows from the trace condition. The converse, i.e. that all these equalities taken together imply $\var$
 is a trace, follows from Lemma 9.4.
\qed

{\bf 9.11.} We will now look at dual-positivity. {\em We shall assume for the rest of the section {\rm 9} that
$E$ and $B$ are $C^*$-algebras and that $Y=Y^*$. Note that $(\rho_n(Y;B))^* = \rho_n(Y;B)$ and $R_n(Y;B)(b)=
(R_n(Y;B)(b^*))^*$ i.e. $(R(Y;B))^* = R(Y;B)$ under these assumptions.}
\\

{\bf 9.12. Proposition.} \\

(i) {\em We have $U(\var^*) = (U(\var))^*$.}\\

(ii) {\em We have $U(\var)^* = U(\var)$ iff $\var = \var^*$.}\\

(iii) {\em $\var\geq 0$ iff $-U(\var)\geq 0$ in the sense of dual positivity in $A(\rho(Y;B))$.\\

Proof.} (i) If $b\in \rho_n(Y;B)$ and $t$ denotes transpose of a matrix, then
\begin{eqnarray*}
(U_n(\var^*))(b)&=& (\var^* \ox {\rm id}_{{\mathfrak M}_n})((b- Y\ox I_n)^{-1})\\
\\
&=& \overline{(\var\ox {\rm id}_{{\mathfrak M}_n})((b-Y\ox I_n)^{-1})^{*t})}\\
\\
&=& \overline{(\var\ox {\rm id}_{{\mathfrak M}_n})((b^*-Y\ox I_n)^{-1})^t}\\
\\
&=&(U_n(\var)(b^*))^*= (U_n(\var))^*(b).
\end{eqnarray*}
 
(ii) follows from (i) and the injectivity of $U$.\\

(iii) We first prove the only if part. Assume $\var \geq 0$ and let $h\in {\mathfrak M}_n, h\geq 0$. By
Lemma 7.13 and the definition of dual positivity, we must check that in the $2n\x 2n$ matrix $(\var \ox
{\rm id}_{{\mathfrak M_{2n}}})((b\oplus b^*- Y\ox I_{2n}-1\ox \g_{n,n}(h))^{-1}) $ the right $n\x n$ corner 
block is positive. Since this block is precisely $(\var \ox
{\rm id}_{{\mathfrak M_n}})((b-Y\ox I_n)^{-1}(1\ox h)(b^* - Y\ox I_n)^{-1})$ the assertion follows from
the assumptions $\var\geq 0$ and $h\geq 0$.

To prove the converse, note that from the proof of the only if part the dual positivity of $-U(\var)$ implies 
$(\var \ox
{\rm id}_{{\mathfrak M_n}})((b-Y\ox I_n)^{-1}(1\ox h)(b^* - Y\ox I_n)^{-1}) \geq 0$ for all
 $h\geq 0, h\in {\mathfrak M}_n$. This in turn implies $\var(\xi \xi^*) \geq 0$ for any $\xi$ in the linear
span  of ${\cal CR}(Y;B)$. Indeed if $\xi=c_1\eta_1+\dots + c_p\eta_p$ where $c_j\in {\mathbb C}$ and $\eta_j$
is  some matrix coefficient of  $(b_j-Y\ox I_{n_j})^{-1}$ then it is easily seen that 
\[
\xi\xi^* = (1\ox k)(b-Y\ox I_n)^{-1}(1\ox h)(b^*-Y\ox I_n)^{-1}(1\ox k^*)
\]
for some $h\geq 0, h\in {\mathfrak M}_n, k\in {\mathfrak M}_{1,n}$ and $n=n_1+\dots +n_p, b=b_1\oplus \dots 
\oplus b_p$.
Hence $\var(\xi\xi^*)= k(\var\ox {\rm id}_{{\mathfrak M}_n})((b-Y\ox I_n)^{-1}(1\ox h)(b^*-Y\ox I_n)^{-1})k^*
\geq 0$.
\qed

{\bf 9.13. Remark.} The dual positivity of $-U(\var)$ is equivalent to the dual positivity of $\overline{U}(\var)$
w.r.t $-\p$, which is then in agreement with $\p$ intertwining the GDQ structures of $E'$ and
$(A(\rho(Y;B)), -\p)$.\\

{\bf 9.14. Remark.} {\em To characterize states in $E$ via their $U$-transform one requires in addition to 
dual-positivity of $-U(\var)$ also $\var(1)=1$, which is equivalent to $\lim_{n\rightarrow \i}nU_1(\var)(n1)
=1$ ($n1\in \rho_1(Y;B)$ for $n\geq \|Y\|$).}\\

{\bf 9.15. Remark.} One situation in free probability where the dual multiplication appears is the definition
of the conjugate variable ${\cal J}(X:B)$ ([11] see also [13], [14]). In the corresponding $W^*$-probability
context $(M, \tau)$ and $B\<X\>\subset M$, with $1\in B$ a von Neumann subalgebra, $\tau$ a trace state and
assuming $B\<X\>$ weakly dense in $M$ let $\var(\cdot) = \tau(\cdot{\cal J}(X:B))$ be the functional defined
by ${\cal J}(X:B)$. Then, if $a\in B\<X\>$ we have $\var(a)=(\tau\ox \tau)(\p_{X:B}a)$ or, denoting by $\#$ the
dual multiplication $\tau \# \tau=\var$. Identifying $L^2(M,\tau)$ with a part of the predual $M_*$ of $M$
 and hence $\tau$ with 1, the same relation would be written in the form $1\# 1 = {\cal J}(X:B)$. Similarly the
higher conjugates ([11]) amount to $(p+1)$ fold dual products $\tau\#\dots\#\tau$ or in the other notation
$1\#\dots\# 1$.

\end{document}